\documentclass[12pt]{amsart}

\usepackage{amsmath,graphics,comment}

\theoremstyle{plain}
\newtheorem*{theorem*}{Theorem}
\newtheorem*{lemma*} {Lemma}
\newtheorem*{corollary*} {Corollary}
\newtheorem*{proposition*} {Proposition}
\newtheorem{theorem}{Theorem}[section]

\newtheorem{corollary}[theorem]{Corollary}
\newtheorem{proposition}[theorem]{Proposition}

\theoremstyle{remark}
\newtheorem*{remark}{Remark}

\theoremstyle{definition}

\def \Z {\mathbf{Z}}

\def \A {\mathbf{A}}
\def \L {\mathbf{L}}

\def\eps{\epsilon}

\def\s{\sigma}

\def\Q{\Bbb{Q}}

\def\Z{\Bbb{Z}}

\def\l{\lambda}
\def\part{\partial}

\def\a{\alpha}

\def\bp{\begin{pmatrix}}

\def\sm{\setminus}
\def\ep{\end{pmatrix}}
\def\bn{\begin{enumerate}}

\def\rk{\mbox{rank}}
\def\en{\end{enumerate}}
\def\ba{\begin{array}}
\def\ea{\end{array}}

\def\L{\Lambda}

\def\s{\sigma}
\def\a{\alpha}

\def\ti{\tilde}
\def\lk{\mbox{lk}}

\def\fr12{\frac{1}{2}}
\def\diag{\mbox{diag}}

\def\ker{\mbox{Ker}}

\def\e{\epsilon}

\begin{document}

\title{Algorithm for finding boundary link Seifert matrices }
\author{Stefan Friedl}
 \address{Rice University, Houston, Texas, 77005-1892}
\email{friedl@math.rice.edu}

\date{\today}
\begin{abstract}
We explain an algorithm for finding a boundary link Seifert matrix for a given multivariable
Alexander polynomial. The algorithm depends on several choices and therefore makes it possible
to find non-equivalent Seifert matrices for a given Alexander polynomial.
\end{abstract}
\maketitle

\section{Introduction}

\subsection{Algebraic statement}
We call $A=(A_{ij})_{i,j=1,\dots,m}$ a (boundary link) Seifert matrix if $A$ is a matrix with
entries $A_{ij}$ which are $(n_i \times n_j)$--matrices over $\Z$ such that $A_{ij}=A_{ji}^t$
for $i \ne j$ and $\det(A_{ii}-A_{ii}^t)=1$ (for more details cf. \cite{L77}, \cite{K87}) Note
that the $n_i$ are necessarily even numbers. Set
\[T:=\diag(\underbrace{t_1,\dots,t_1}_{n_1},
\dots,\underbrace{t_m,\dots,t_m}_{n_m}), \] then define the Alexander polynomial of $A$ to be
\[\Delta(A):=\det(T)^{-\frac{1}{2}}\det(TA-A^t) \in \L_m:=\Z[t_1^{\pm 1},\dots,t_m^{\pm 1}].\]
This polynomial has the following well-known properties
which can easily be verified from the definitions.
\[ \ba{rcl} \Delta(A)(1,\dots,1)&=&1, \\
 \Delta(A)(t_1,\dots,t_m)&=&\Delta(A)(t_1^{-1},\dots,t_m^{-1}). \ea \]

Now assume that $\Delta$ is a polynomial with the above properties. The goal of this paper is to
give an algorithm how for finding a Seifert matrix $A$ in terms of the coefficients of $\Delta$
such that $\Delta(A)=\Delta$. In the case $m=1$, i.e. the case of Seifert matrices for knots, an
algorithm has been found by Seifert (cf. \cite{S34}, \cite{BZ85}).

\subsection{Topological motivation}

We quickly recall how boundary link Seifert matrices appear
in link theory.
An $m$-link $L=L_1\cup\dots \cup L_m \subset S^{4q+3}$
is a smooth embedding of
$m$ disjoint oriented $(4q+1)$-spheres. A boundary link is a link
which has $m$ disjoint Seifert manifolds, i.e. there exist
$m$ disjoint oriented $(4q+2)$-submanifolds $F_1,\dots,F_m \subset S^{4q+3}$
such that $\partial(F_i)=L_i, i=1,\dots,m$.
One of the main tools
for studying boundary links is the Seifert form
\[ \ba{rcl}
H_{2q+1}(F) \times H_{2q+1}(F) &\to & \Z \\
(a,b) &\mapsto &\lk(a,b_+),
\ea \]
where $b_+$ means that we push a representative of $b$
into $S^{4q+3}\sm F$ along the positive normal direction of $F$.
More precisely, we can find an orientation preserving embedding $\iota:F\times [-1,1] \to S^{4q+3}$
and we define $a_+=\iota(a,+1)$ and $a_-=\iota(a,-1)$.

Now pick bases $l_{i,1},\dots,l_{i,n_i}$
for $H_{2q+1}(F_i), i=1,\dots,m$,
then
\[ l_{1,1},\dots,l_{1,n_1},\dots,l_{m,1},\dots,l_{m,n_m}\]
form a basis for $H_{2q+1}(F)=H_{2q+1}(F_1)\oplus \dots \oplus H_{2q+1}(F_m)$.
Representing the Seifert form with respect to this basis
we get a boundary link Seifert matrix (cf. \cite{L77}, \cite[p. 670]{K87}).

We also need the notion of an $F_m$--link, this is a link with a map $\pi_1(S^{4q+3}\sm L)\to
F_m$, where $F_m$ denotes the free group on $m$ generators, which sends meridians to conjugates
of the generators of $F_m$. A Thom argument shows that there is a one-to-one correspondence
between isotopy classes of $F_m$--links and isotopy classes of boundary links with Seifert
manifolds. It turns out that it is easier to study $F_m$--links, for example the addition of
$F_m$--links is well-defined for $q\geq 1$. Boundary links and $F_m$--links are the best
understood links, they have been studied thoroughly and many of the classifying results for
higher dimensional knots can be done similarly in the context of such links (cf. \cite{L77},
\cite{K87}, \cite{D86}).

If $L$ is a boundary link with $m$ components then denote by $\ti{X}$ the universal abelian
cover of $S^{4q+3}\sm L$, i.e. the cover induced by
$\pi_1(S^{4q+3}\sm L)\to H_1(S^{4q+3}\sm L)=\Z^m$.
Note that $H_*(\ti{X})$ has a natural $\Z[\Z^m]=\L_m$-module
structure.

\begin{proposition}
Let $L \subset S^{4q+3}$ be a   boundary link with $m$ components, and $A$ a Seifert matrix of
size $(n_1,\dots,n_m)$ for a Seifert manifold $F=F_1\cup \dots\cup F_m$. Then there exists a
short exact sequence
\[ 0\to \L_m^n/(AT-A^t)\L_m^n\to  H_{2q+1}(\ti{X}) \to  P\to 0,\]
where $n=\sum_{i=1}^m n_i$ and  $P$ is some torsion free $\L_m$-module. Furthermore $P=0$ if
$q>0$.
\end{proposition}

We will give a quick outline of the proof which follows well-known arguments in the knot case
(cf. \cite{L66}, \cite{R90}).

\begin{proof}
 Let
$Y=S^{4q+3}\sm F$.
We can view $\ti{X}$ as the result of gluing $\Z^m$ copies of $Y$ together along $\Z^m$ copies of
$F_1,\dots,F_m$.
Consider the resulting Mayer-Vietoris sequence
\[
\ba{rcl}
 \dots \to H_{i}(F)\otimes \L_m& \to& H_{i}(Y)\otimes \L_m \to H_{i}(\ti{X}) \to \dots \\
       a_j \otimes p &\mapsto & (a_{j,+}t_j-a_{j,-})\otimes p,
\ea \]
where $a_j \in H_i(F_j)$.
Note that for $i\in \{1,\dots,4q+1\}$
we have
$H_i(Y)\cong H^{4q+2-i}(F)\cong H_i(F,L)\cong H_i(F)$ by Alexander duality, Poincar\'e duality
and a long exact sequence argument.
 Pick a basis for $H_i(F)$ which gives $A$ as a Seifert matrix for $L$, then give
$H_i(Y)$ the corresponding basis.
An argument as in Rolfsen (cf. \cite{R90}) shows that the map $H_{2q+1}(F)\otimes \L_m
\to H_{2q+1}(Y)\otimes \L_m$
is given by $v\mapsto (AT-A^t)v$.

If $q=0$ then the sequence becomes
\[ \dots \to H_{1}(F)\otimes \L_m \to H_{1}(F)\otimes \L_m \to H_{1}(\ti{X}) \to
 H_{0}(F)\otimes \L_m \to H_{0}(Y)\otimes \L_m \to \]
It is clear that $\ker\{H_{0}(F)\otimes \L_m \to H_{0}(Y)\otimes \L_m\}$ is $\L_m$-torsion free
since $H_0(Y)$ is $\L_m$-torsion free (cf. also \cite{S81b}).

Now consider the case $q>0$. We are done once we show that $H_{2q}(F)\otimes \L_m \to
H_{2q}(Y)\otimes \L_m\cong H_{2q}(F)\otimes \L_m$ is injective. Picking a basis for $H_{2q}(F)$
and giving $H_{2q}(Y)$ the corresponding basis, then we can represent this map by a matrix
$B(t_1,\dots,t_m)$. We will prove that $B(1,\dots,1)$ is in fact the identity matrix, in
particular $\det(B)\ne 0$. This concludes the proof of the proposition. Note that $B(1,\dots,1)$
represents the map $H_{2q}(F) \to H_{2q}(Y)$ given by $a\mapsto a_+-a_-$. Recall that the
isomorphism $f:H_i(Y) \to H^{2q+2-i}(F)$ is induced by the linking pairing, in particular for
$\s\in C_{2q+2-i}(F)$
\[ f(a_+-a_-)(\s)=\lk(a_+-a_-,\s)=(a\times [-1,1])\cdot \s=a \cdot \s.\]
Thus under the Poincar\'e duality map $f(a_+-a_-)$ gets sent to $a$.
\end{proof}

From the theory of fitting ideals for presentation matrices (cf. \cite{S81b}) it follows that
$\det(AT-A^t)$ is a well-defined invariant for a boundary link $L$ up to multiplication by a
unit in $\L_m$. It is easy to see that $\det(T)^{-\frac{1}{2}}\det(AT-A^t)$ is a well-defined
invariant for boundary links, it is called the Alexander polynomial of $L$.

Gutierrez \cite[p. 34]{G74} showed that any polynomial
$\Delta(t_1,\dots,t_m)$ with the properties
\[ \ba{rcl} \Delta(1,\dots,1)&=&1 \\
 \Delta(t_1,\dots,t_m)&=&\Delta(t_1^{-1},\dots,t_m^{-1}) \ea \]
is the Alexander polynomial of a boundary link in dimension 1, in particular there exists a
boundary link Seifert matrix $A$ with $\Delta(A)=\Delta$. But it is   difficult to find an
explicit boundary link Seifert matrix, which would be important to compute further invariants.

\begin{remark}
Farber \cite{F92} and Garoufalidis and Levine \cite{GL02} defined non-commutative invariants
for boundary links which can be viewed as generalizations of the Alexander polynomial of a knot.
Farber also proves a realization theorem.
\end{remark}

\subsection{$S$-equivalence class of Seifert matrices}
In the following we will call a matrix $P$ block diagonal if it commutes with $T$, equivalently
if $P=P_1\oplus \dots \oplus P_m$ where $P_i$ is a $(n_i\times n_i)$-matrix.

The $S$-equivalence of Seifert matrices is the equivalence relation  generated by the following
two equivalences (for more details cf. \cite{L77}, \cite{K85})
\bn
\item $A \sim PAP^t$ where $P$ is a block diagonal matrix over $\Z$ with
$\det(P)=1$.
\item $A$ is equivalent to any row or column enlargement or reduction
of $A$.
\en

\begin{proposition}[\cite{L77}, \cite{K85}] \label{propko}
Any two Seifert matrices for an $F_m$--link are $S$-equivalent.
Furthermore any Seifert matrix is the Seifert matrix of an $F_m$--link.
\end{proposition}

There exists a similar but more complicated proposition for boundary links
(cf. \cite{K85}). It turns out that Seifert matrices for boundary links are related
by $S$-equivalence and an action by (cf. also \cite{K87})
\[ A_m:=\{\varphi:F_m\to F_m | \varphi(x_i)=\l_i x_i\l_i^{-1} \mbox { for some }\l_i\in F_m\}/\mbox{inner
automorphism}.\] The groups $A_1,A_2$ are trivial \cite{K84}, it follows that boundary link
matrices with 2 components which are are related by $S$-equivalence and an action by $A_m$ are
in fact $S$-equivalent.

It is easy to see that if $A_1,A_2$ are $S$-equivalent, then $\Delta(A_1)=\Delta(A_2)$, this
shows again that the Alexander polynomial is an invariant for any $F_m$--link.

We call a Seifert matrix irreducible if no row or column reductions are possible.

\begin{proposition} \label{propprops}
\bn
\item
A Seifert matrix of size $(n_1,\dots,n_m)$ is irreducible if and only if
\[\rk(A_{i1} \dots A_{im})=n_i,
  \quad \rk\left(\ba{c} A_{1i}\\ \vdots \\A_{mi}\ea \right)=n_i\]
for all $i=1,\dots,m$. Put differently, a Seifert matrix is irreducible
if and only if the block columns and block rows have maximal rank.
\item If $A_1,A_2$ are $S$-equivalent minimal Seifert matrices then
$A_1=PA_2P^t$ where $P$ is a block diagonal matrix over $\Q$ with
$\det(P)\ne 0$.
\en
\end{proposition}

We will use this proposition to show that certain Seifert matrices are not $S$-equivalent.

The statement of the proposition is well-known in the case $m=1$ (cf. \cite{T73}). The first
part of the proposition is fairly straight forward to show, whereas the second part is more
difficult to prove. Using ideas of Farber \cite{F91} one can rewrite the proof of Trotter in the
general case, but this requires many details, which we will omit here.


\section{Statement of results}
\subsection{Algebra}
For $v_1,\dots,v_l \in \Z$ and $\e_2,\dots,\e_l \in \{-1,+1\}$
define matrices $B_i:=B_i(v_1,\dots,v_i,\eps_2,\dots,\eps_i)$
inductively as follows.
\[
  B_1 :=\bp v_1 & 0 \\-1 & 1 \ep \quad
B_i := \bp &&&&&&& v_i & 0 \\
&&&&&&& 0 & 1 \\&&&&&&& v_i & 1 \\&&&&B_{i-1}&&& 0 & 1 \\
&&&&&&& \vdots \\
&&&&&&& v_i & 1 \\&&&&&&& 0 & 1 \\
v_i &0&v_i&0&\dots &v_i &0&v_i&z_n \\
0&1&1&1 &&1&1&1-z_n&1 \ep,
\]
where $z_i:=\frac{1}{2}(1+\e_i)$.
Furthermore let
\[ Y_l:=\diag(y_1,y_1,y_2,y_2,\dots,y_l,y_l). \]

\begin{proposition} \label{mainprop}
Set $v_{l+1}=0$, then
\[ \det(Y_lB_l-Y_l^{-1}B_l^t)=1-2v_1+ \sum_{j=1}^l (v_j-v_{j+1})
\left(y_1^2\prod_{i=2}^j y_i^{2\e_i}+y_1^{-2}\prod_{i=2}^j y_i^{-2\e_i}\right). \]
\end{proposition}

The proof will be given in Section \ref{sectionproof}.

\subsection{Explanation of the algorithm}
Let $\Delta \in \Z[t_1,\dots,t_m]$ be a polynomial with the
following properties
\[ \ba{rcl} \Delta(1,\dots,1)&=&1, \\
 \Delta(t_1,\dots,t_m)&=&\Delta(t_1^{-1},\dots,t_m^{-1}). \ea \]
Then using the usual multiindex notation we can uniquely write
\[ \Delta(t_1,\dots,t_m)=\sum_{\a \in \Z^m}c_{\a}(t^{\a}+t^{-\a})
+1-\sum_{\a \in \Z^m}2c_{\a}, \quad c_{\a}\in \Z,\] where $c_{\a}=0$ for all but finitely many
$\a$ and $c_{(0,\dots,0)}=0$.

Denote the $\a$ with $c_{\a}\ne 0$ by $\a_1,\dots,\a_r$.
Pick a map $p:\{0,\dots,l\}\to \Z^m$ with the following properties.
\bn
\item $p(0)=(0,\dots,0)$,
\item $|p(t)-p(t-1)|=1$ for all $t=1,\dots,l$,
\item for each $i=1,\dots,r$
there exists a $t_i\in \{1,\dots,l\}$ such that $p(t_i)=\a_i$.
\en
It is easy to see that such a map always exists. Denote the $i^{th}$ unit vector in $\Z^m$ by
$e_i$, the second condition says that $p(t)=p(t-1)+\e_te_{s_t}$ for unique $\e_t \in \{-1,+1\},
s_t\in \{1,\dots,m\}$.

Now define $w_{t_i}=c_{p(t_i)}=c_{\a_i}$ for $i=1,\dots,r$ and $w_j=0$ otherwise. Let
$v_i:=\sum_{j=i}^{l}w_j, j=1,\dots,l$. From Proposition \ref{mainprop} it follows now
immediately that for $B=B(v_1,\dots,v_l,\e_2,\dots,\e_l)$ and $Y:=\diag(y_{s_1},\dots,y_{s_l})$
we get
\[ \ba{rcl}
\det(YB-Y^{-1}B^t)&=& 1-2v_1+\sum_{j=1}^l (v_j-v_{j+1}) \left(y_{s_1}^2\prod_{i=2}^j
y_{s_i}^{2\e_i} +y_{s_1}^{-2}\prod_{i=2}^j y_{s_i}^{-2\e_i}\right).\ea \] Using multiindex
notation $y=(y_1,\dots,y_m)$ we can rewrite this as
\[\sum_{j=1}^l w_j(y^{p(j)}+y^{-p(j)})
+1-\sum_{j=1}^l 2w_j = \sum_{\a \in \Z^m}c_{\a}(y^{2\a}+y^{-2\a}) +1-\sum_{\a \in \Z^m}2c_{\a},
\] in particular for $\ti{T}:=\diag(t_{s_1},\dots,t_{s_l})$ we get
\[ \det(\ti{T})^{\frac{1}{2}}\det(\ti{T}B-B^t)
=\sum_{\a \in \Z^m}c_{\a}(t^{\a}+t^{-\a}) +1-\sum_{\a \in \Z^m}2c_{\a}=\Delta. \] We can find a
permutation matrix $P$ such that
\[ P\ti{T}P^{-1}=\diag(\underbrace{t_1,\dots,t_1}_{n_1},
\dots,\underbrace{t_m,\dots,t_m}_{n_m})=:T \]
for some $n_1,\dots,n_m$.
In fact we can and will assume that $P$ is of form
\[ P(v_{1,1},v_{1,2},v_{2,1},v_{2,2},\dots,v_{l,1},v_{l,2})
=P(v_{\s(1),1},v_{\s(1),2},v_{\s(2),1},v_{\s(2),2},
\dots,v_{\s(l),1},v_{\s(l),2})
\]
for some permutation $\s \in S_l$, i.e. $P$ permutes pairs
of coordinates. Note that $P^t=P^{-1}$ and $\det(P)=1$.

\begin{theorem}
The matrix $A=PBP^{-1}$ is a boundary link Seifert matrix
of size $(n_1,\dots,n_m)$
and $\Delta(A)=\Delta$.
\end{theorem}

\begin{proof}
Note that $B-B^t$ and hence $A-A^t$ is a block sum of $2 \times 2$ matrices
of the form $\bp 0 & \pm 1 \\ \mp 1 &0 \ep$, in particular
$A$ is a Seifert matrix of size $(n_1,\dots,n_m)$,
furthermore
\[\ba{rcl} \Delta(A)
&=&\det({T})^{-\frac{1}{2}}\det(TA-A^t)
=\det({T})^{-\frac{1}{2}}\det(PTP^{-1}PAP^{-1}-PA^tP^{-1})=\\
&=&\det(\ti{T})^{-\frac{1}{2}}\det(\ti{T}B-B^t)=\Delta. \ea \]
\end{proof}

Using Proposition \ref{propko} we get the following corollary.

\begin{corollary}
Any  $\Delta$ with $\Delta(1,\dots,1)=1$
and $\Delta(t_1^{-1},\dots,t_m^{-1})=\Delta(t_1,\dots,t_m)$
is the Alexander polynomial of a boundary link.
\end{corollary}

It is clear that $A$ depends on the map $p$, for example $A$ is a $(2l \times 2l)$--matrix, i.e.
$p$ determines the size of $A$. We will see in the next section that different paths can in fact
give non $S$-equivalent matrices.

\subsection{Example}
\subsubsection{Minimality of matrices}
Let $\Delta=c_{1,0}(t_1+t_1^{-1})+c_{1,1}(t_1t_2+t_1^{-1}t_2^{-1})+c_{0,1}(t_2^2+t_2^{-2})-17$,
then $\a_1=(1,0), \a_2=(1,1), \a_3=(0,1)$.
The map $p(0):=(0,0), p(1):=(1,0), p(2):=(1,1), p(3):=(0,1)$
satisfies the conditions on $p$. In this case
\[ \ba{rclrclrcl}
t_1&=&1,& t_2&=&2,&t_3&=&3,\\
s_1&=&1,& s_2&=&2,&s_3&=&1,\\
\e_1&=&1,& \e_2&=&1,&\e_3&=&-1,\\
w_1&=&c_{1,0},& w_2&=&c_{1,1},&w_3&=&c_{0,1},\\
v_1&=&c_{1,0}+c_{1,1}+c_{0,1},& v_2&=&c_{1,1}+c_{0,1},&v_3&=&c_{0,1}.
\ea
\]
Then
\[ B=\bp v_1&0\hspace{0.2cm}&v_2&0\hspace{0.2cm}&v_3&0 \\
-1&1\hspace{0.2cm}&0&1\hspace{0.2cm}&0&1
\\[0.2cm]v_2&0\hspace{0.2cm}&v_2&1\hspace{0.2cm}&v_3&1
\\0&1\hspace{0.2cm}&0&1\hspace{0.2cm}&0&1
\\[0.2cm]v_3&0\hspace{0.2cm}&v_3&0\hspace{0.2cm}&v_3&0
\\0&1\hspace{0.2cm}&1&1\hspace{0.2cm}&1&1
\ep \mbox{ and }
A=\bp v_1&0\hspace{0.2cm}&v_3&0\hspace{0.2cm}&v_2&0 \\
-1&1\hspace{0.2cm}&0&1\hspace{0.2cm}&0&1
\\[0.2cm]v_3&0\hspace{0.2cm}&v_3&0\hspace{0.2cm}&v_3&0
\\0&1\hspace{0.2cm}&1&1\hspace{0.2cm}&1&1
\\[0.2cm]v_2&0\hspace{0.2cm}&v_3&1\hspace{0.2cm}&v_2&1
\\0&1\hspace{0.2cm}&0&1\hspace{0.2cm}&0&1
\ep, \] where we chose $ \s= \left(\ba{ccc} 1 &2&3 \\ 1&3&2\ea \right)$.

Using Proposition \ref{propprops} it is easy to see that $A$ forms an irreducible Seifert matrix
of size $(2,1)$.

Consider
\[ A=\bp w_1+w_3 &0\hspace{0.2cm} &-w_3&0 \\ -1&1\hspace{0.2cm}&0&1 \\[0.2cm]
-w_3&0\hspace{0.2cm}&w_2+w_3&0 \\ 0&1\hspace{0.2cm}&-1&1 \ep \]
then
\[ \Delta(A)=w_1(t_1+t_1^{-1})+w_2(t_2+t_2^{-1})+w_3(t_1t_2+t_1^{-1}t_2^{-1})+1-2(w_1+w_2+w_3).\]
This shows that the algorithm does in general not produce a Seifert matrix of minimal size
for a given Alexander polynomial.

\subsubsection{Uniqueness of result}
A straight forward argument shows that for a knot Alexander polynomial $\Delta(t)$
different choices of maps $p$ will produce $S$-equivalent matrices. This is no longer true in
the case $m>1$.

Consider $\Delta=w(t_1t_2+t_1^{-1}t_2^{-1})+1-2w, w\ne 0$. If we take maps $p_1,p_2$ with
$p_1(0)=(0,0), p_1(1)=(1,0)$ and $p_1(2)=(1,1)$ and $p_2(0)=(0,0), p_2(1)=(0,1)$ and
$p_2(2)=(1,1)$ then applying the algorithm we will get identical matrices $B$ but we have to use
different permutations:
\[ \s_1 =\bp 1 & 2 \\ 1&2 \ep, \quad \s_1 =\bp 1 & 2 \\ 2&1 \ep. \]
We get Seifert matrices\[
 A_1=\bp w&1\hspace{0.2cm}&w&0 \\ 0&1\hspace{0.2cm}&1&1\\[0.2cm]
w&1\hspace{0.2cm}&0&1\\ 0&1\hspace{0.2cm}&0&1 \ep
\mbox{ and }
A_2=\bp 0&1\hspace{0.2cm}&w&1 \\ 0&1\hspace{0.2cm}&0&1\\[0.2cm]
w&0\hspace{0.2cm}&w&1\\ 1&1\hspace{0.2cm}&0&1 \ep.
\]
Both matrices are minimal, but not block congruent, since $\det(A_{1,11})=w, \det(A_{2,11})=0$
Hence by Proposition \ref{propprops} $A_1$ and $A_2$ are not $S$-equivalent.

Recall that any boundary link Seifert matrix corresponds to an $F_m$--link, we therefore can
construct non-isotopic $F_m$--links with identical Alexander polynomials. I do not know whether
the matrices are $S_m$-equivalent, in particular whether the corresponding boundary links are
isotopic.

Using signature invariants one can show that these matrices are in fact not even matrix
cobordant (for a definition cf. \cite{K87}), i.e. one can show that the corresponding
$F_m$--links are in fact not even $F_m$--cobordant.

\section{Proof of Proposition \ref{mainprop}} \label{sectionproof}

\subsection{Proof of a special case of Proposition \ref{mainprop}}
In this section we will consider the case $\eps_2=\dots=\eps_l=1$. We have to show that
\[ \det(Y_lB_l-Y_l^{-1}B_l^t)=1-2v_1+ \sum_{j=1}^l (v_j-v_{j+1})
\left(\prod_{i=1}^j y_i^2+\prod_{i=1}^j y_i^{-2}\right). \] We will show how to compute the
determinant, but we will give the matrices only for the case $l=4$ to simplify the notation.

Consider $Y_lB_l-Y_l^{-1}B_l^t$:
{\footnotesize
\[
\bp
v_1(y_1-y_1^{-1})\hspace{-0.3cm}&y_1^{-1}
&v_2(y_1-y_1^{-1})\hspace{-0.3cm}&0&v_3(y_1-y_1^{-1})\hspace{-0.3cm}&0&v_4(y_1-y_1^{-1})&0 \\
-y_1\hspace{-0.3cm}&y_1-y_1^{-1}&0\hspace{-0.3cm}&y_1-y_1^{-1}&0\hspace{-0.3cm}&y_1-y_1^{-1}&0&y_1-y_1^{-1}
\\[0.2cm]
v_2(y_2-y_2^{-1})\hspace{-0.3cm}&0&v_2(y_2-y_2^{-1})\hspace{-0.3cm}&y_2&v_3(y_2-y_2^{-1})\hspace{-0.3cm}&y_2-y_2^{
-1}
&v_4(y_2-y_2^{-1})&y_2-y_2^{-1} \\
0\hspace{-0.3cm}&y_2-y_2^{-1}&-y_2^{-1}\hspace{-0.3cm}&y_2-y_2^{-1}&0\hspace{-0.3cm}&y_2-y_2^{-1}&0&y_2-y_2^{-1}
\\[0.2cm]
v_3(y_3-y_3^{-1})\hspace{-0.3cm}&0&v_3(y_3-y_3^{-1})\hspace{-0.3cm}&0&v_3(y_3-y_3^{-1})\hspace{-0.3cm}&y_3&v_4(y_3
-y_3^{-1})&y_3-y_3^{-1} \\
0\hspace{-0.3cm}&y_3-y_3^{-1}&y_3-y_3^{-1}\hspace{-0.3cm}&y_3-y_3^{-1}&-y_3^{-1}\hspace{-0.3cm}&y_3-y_3^{-1}&0&y_3
-y_3^{-1} \\[0.2cm]
v_4(y_4-y_4^{-1})\hspace{-0.3cm}&0&v_4(y_4-y_4^{-1})\hspace{-0.3cm}&0&v_4(y_4-y_4^{-1})\hspace{-0.3cm}&0&v_4(y_4-y
_4^{-1})&y_4 \\
0\hspace{-0.3cm}&y_4-y_4^{-1}&y_4-y_4^{-1}\hspace{-0.3cm}&y_4-y_4^{-1}&y_4-y_4^{-1}\hspace{-0.3cm}&y_4-y_4^{-1}&-y
_4^{-1}&y_4-y_4^{-1} \\
\ep.
\]
} We will first simplify the matrix to make the computation of the determinant easier. For
$i=2,\dots,l$ multiply the second row by $\frac{y_i-y_i^{-1}}{y_1-y_1^{-1}}$ and subtract the
result from the $2i$--th row, we get {\footnotesize
\[ \bp
v_1(y_1-y_1^{-1})\hspace{-0.3cm}&y_1^{-1}
&v_2(y_1-y_1^{-1})\hspace{-0.3cm}&0&v_3(y_1-y_1^{-1})\hspace{-0.3cm}&0&v_4(y_1-y_1^{-1})&0 \\
-y_1\hspace{-0.3cm}&y_1-y_1^{-1}&0\hspace{-0.3cm}&y_1-y_1^{-1}&0\hspace{-0.3cm}&y_1-y_1^{-1}&0&y_1-y_1^{-1}
\\[0.2cm]
v_2(y_2-y_2^{-1})\hspace{-0.3cm}&0&v_2(y_2-y_2^{-1})\hspace{-0.3cm}&y_2&v_3(y_2-y_2^{-1})\hspace{-0.3cm}&y_2-y_2^{
-1}&v_4(y_2-y_2^{-1})&y_2-y_2^{-1} \\
y_1\frac{y_2-y_2^{-1}}{y_1-y_1^{-1}}\hspace{-0.3cm} &0&-y_2^{-1}\hspace{-0.3cm}&0&0\hspace{-0.3cm}&0&0&0 \\[0.2cm]
v_3(y_3-y_3^{-1})\hspace{-0.3cm}&0&v_3(y_3-y_3^{-1})\hspace{-0.3cm}&0&v_3(y_3-y_3^{-1})\hspace{-0.3cm}&y_3&v_4(y_3
-y_3^{-1})&y_3-y_3^{-1} \\
y_1\frac{y_3-y_3^{-1}}{y_1-y_1^{-1}}\hspace{-0.3cm}&0

&y_3-y_3^{-1}\hspace{-0.3cm}&0&-y_3^{-1}\hspace{-0.3cm}&0&0&0 \\[0.2cm]
v_4(y_4-y_4^{-1})\hspace{-0.3cm}&0&v_4(y_4-y_4^{-1})\hspace{-0.3cm}&0&v_4(y_4-y_4^{-1})&0&v_4(y_4-y_4^{-1})\hspace
{-0.3cm}&y_4 \\y_1\frac{y_4-y_4^{-1}}{y_1-y_1^{-1}}\hspace{-0.3cm}&0&y_4-y_4^{-1}\hspace{-0.3cm}&0
&y_4-y_4^{-1}\hspace{-0.3cm}&0&-y_4^{-1}&0 \\
\ep.
\]
} For $i=1,\dots,l-1$ subtract the $(2i+1)$--st column from the $(2i-1)$--st column and for
$i=l-1,\dots,1$ subtract the $2i$--th column from the $(2i+2)$--nd column, we get {\footnotesize
\[
\bp w_1(y_1-y_1^{-1})\hspace{-0.3cm}&y_1^{-1}&w_2(y_1-y_1^{-1})\hspace{-0.3cm}&-y_1^{-1}&
w_3(y_1-y_1^{-1})&0&w_4(y_1-y_1^{-1})&0 \\
-y_1\hspace{-0.3cm}&y_1-y_1^{-1}&0\hspace{-0.3cm}&0&0\hspace{-0.3cm}&0&0&0\\[0.2cm]
0\hspace{-0.3cm}&0&w_2(y_2-y_2^{-1})\hspace{-0.3cm}&y_2&w_3(y_2-y_2^{-1})\hspace{-0.3cm}&-y_2^{-1}&w_4(y_2-y_2^{-1
})&0\\
y_1\frac{y_2-y_2^{-1}}{y_1-y_1^{-1}}+y_2^{-1}\hspace{-0.3cm}&0&-y_2^{-1}\hspace{-0.3cm}&0&0\hspace{-0.3cm}&0&0&0\\[0.2cm]
0\hspace{-0.3cm}&0&0\hspace{-0.3cm}&0&w_3(y_3-y_3^{-1})\hspace{-0.3cm}
&y_3&w_4(y_3-y_3^{-1})\hspace{-0.3cm}&-y_3^{-1}\\
y_1\frac{y_3-y_3^{-1}}{y_1-y_1^{-1}}-(y_3-y_3^{-1})\hspace{-0.3cm}&0&y_3\hspace{-0.3cm}&0&-y_3^{-1}\hspace{-0.3cm}
&0&0&0\\[0.2cm]
0\hspace{-0.3cm}&0&0\hspace{-0.3cm}&0&0\hspace{-0.3cm}&0&w_4(y_4-y_4^{-1})&y_4\\
y_1\frac{y_4-y_4^{-1}}{y_1-y_1^{-1}}-(y_4-y_4^{-1})\hspace{-0.3cm}&0&0\hspace{-0.3cm}&0&y_4\hspace{-0.3cm}&0&-y_4^
{-1}&0 \ep,
\]
} where $w_i:=v_i-v_{i+1}, i=1,\dots,l-1$, recall that $v_{l+1}=0$ hence $w_l:=v_l$. For
$i=2,\dots,l$ multiply the $(2i-1)$--st row by $\frac{y_{i-1}-y_{i-1}^{-1}}{y_i-y_i^{-1}}$ and
subtract the result from the $(2i-3)$--rd row, furthermore for $i=2,\dots,l-1$ multiply the
$2i$--th row by $y_iy_{i+1}$ and subtract the result from the $(2i+2)$--nd row. An induction
argument shows that the result is a matrix $D_l$ which is inductively defined as follows.
\[ \ba{rcl}
D_1&=&\bp w_1(y_1-y_1^{-1})&y_1^{-1}\\-y_1&y_1-y_1^{-1}\ep \\[0.5cm]
D_2&=&\bp D_1 &&0&\frac{y_1^{-1}y_2^{-1}-y_1y_2}{y_2-y_2^{-1}} \\
&&0&0\\
0&0&w_2(y_2-y_2^{-1})&y_2\\
\frac{-y_1^{-1}y_2^{-1}+y_1y_2}{y_1-y_1^{-1}}&0&-y_2^{-1}&0 \ep
\ea
\]
and for $n=3,\dots,l$
\[ D_n =\bp
&&&&& 0&0 \\
&&&&& \vdots &\vdots \\
&&&&& 0&0 \\
&&D_{n-1}&&& 0&\frac{y_{n-1}^{-1}(y_{n-2}-y_{n-2}^{-1})}{y_{n-1}-y_{n-1}^{-1}}\\
&&&&& 0&0 \\
&&&&& 0&\frac{y_{n-1}^{-1}y_n^{-1}-y_{n-1}y_n}{y_{n}-y_{n}^{-1}}\\
&&&&& 0&0\\
0&0&0&\dots&0&w_n(y_n-y_n^{-1})&y_n \\
\frac{-y_1^{-1}y_n^{-1}+y_1y_2^2\cdot \dots \cdot y_{n-1}^2y_n}
{y_1-y_1^{-1}}&0&0&\dots&0&-y_n^{-1}&0 \ep.
\]
Note that $\det(D_l)=\det(Y_lB_l-Y_l^{-1}B_l^t)$, we will now compute $\det(D_l)$. For
$n=2,\dots,l$ we denote by $D_n'$ respectively $D_n''$ the matrix obtained from $D_n$ by
deleting the first column and the $(2n-3)$--rd respectively $(2n-1)$--st row. Define
\[ {\det}_n:=\det(D_n), \quad {\det}_n ':=\det(D_n '), \quad {\det}_n'':=\det(D_n''). \]
Using the last row to compute $\det(D_n)$ we get
\[
\hspace{-0.5cm}
\ba{rcl}{\det}_n\hspace{-0.35cm}&=&\hspace{-0.35cm}{\det}_{n\hspace{-0.05cm}-\hspace{-0.05cm}1}
\hspace{-0.1cm}-\hspace{-0.1cm}w_n\hspace{-0.05cm}(y_{n}
\hspace{-0.1cm}-\hspace{-0.1cm}y_{n}^{\hspace{-0.1cm}-1}\hspace{-0.05cm})
\hspace{-0.1cm}\frac{-y_1^ {\hspace{-0.1cm}-1} y_n^{\hspace{-0.1cm}-1}+y_1y_2^2\cdot \dots \cdot
y_{n-1}^2y_n}{y_1-y_1^{\hspace{-0.1cm}-1}} \hspace{-0.15cm}\left(
\frac{y_{n\hspace{-0.05cm}-\hspace{-0.05cm}1}^{\hspace{-0.1cm}-1}
(y_{n\hspace{-0.05cm}-\hspace{-0.05cm}2}-y_{n\hspace{-0.05cm}
-\hspace{-0.05cm}2}^{\hspace{-0.1cm}-1})}
{y_{n\hspace{-0.05cm}-\hspace{-0.05cm}1}-y_{n\hspace{-0.05cm}-
\hspace{-0.05cm}1}^{\hspace{-0.1cm}-1}} {\det}_{n\hspace{-0.05cm}-\hspace{-0.05cm}1}
'\hspace{-0.1cm}+\hspace{-0.1cm}\frac{y_{n\hspace{-0.05cm}-
\hspace{-0.05cm}1}^{-1}y_n^{\hspace{-0.1cm}-1}-y_{n-1}y_n}
{y_n-y_n^{\hspace{-0.1cm}-1}}{\det}_{n\hspace{-0.05cm}-\hspace{-0.05cm}1} ''
\hspace{-0.1cm}\right).
\ea
\]

We make the following easy observations:
\[ \ba{rcl}
{\det}_n '&=&{\det}_{n-1}'',\\
{\det}_n''&=&-y_n^{-1}\left( \frac{y_{n-1}^{-1}(y_{n-2}-y_{n-2}^{-1})}{y_{n-1}-y_{n-1}^{-1}}
{\det}_{n-1} '+\frac{y_{n-1}^{-1}y_n^{-1}-y_{n-1}y_n}{y_n-y_n^{-1}}{\det}_{n-1} '' \right).
\ea \]
It follows that
\[ {\det}_n={\det}_{n-1}-w_n(y_{n}-y_{n}^{-1})
\frac{-y_1^{-1}y_n^{-1}+y_1y_2^2\cdot \dots \cdot y_{n-1}^{2}y_n}{y_1-y_1^{-1}} y_n{\det}_n''.
\]
Recall that we have to show that
\[ \ba{rcl}\det_n&=&1-2v_1+ \sum_{j=1}^l (v_j-v_{j+1})
\left(\prod_{i=1}^j y_i^2+\prod_{i=1}^j y_i^{-2}\right)\\
&=&1-2\sum_{j=1}^lw_j+ \sum_{j=1}^l w_j \left(\prod_{i=1}^j y_i^2+\prod_{i=1}^j y_i^{-2}\right).
\ea \] The proof of the special case of Proposition \ref{mainprop} is complete once we show that
\[ \ba{rcl}
{\det}_1\hspace{-0.1cm}&=\hspace{-0.1cm}&w_1(y_1^2+y_1^{-2})+1-2w_1\\
\frac{-y_1^{-1}y_n^{-1}+y_1y_2^2\cdot \dots \cdot y_{n-1}^{2}y_n}{y_1-y_1^{-1}}
y_n(y_n\hspace{-0.1cm}-\hspace{-0.1cm}y_n^{-1}){\det}_n''\hspace{-0.1cm}&=\hspace{-0.1cm}&y_1^2\cdot
\dots \cdot y_n^2+ y_1^{-2}\cdot \dots \cdot y_n^{-2}\hspace{-0.1cm}-\hspace{-0.1cm}2 \mbox{ for
} n\hspace{-0.1cm}=\hspace{-0.1cm}2,\dots,l.
\ea
\]
The first equality follows from a simple computation.
We now prove the second equality by induction on $n$.
For $n=1,2$ this follows again from a direct computation.
Now assume that the statement is true for all $k<n$, then
using the above results we get
\[
\ba{rl}
&
\frac{-y_1^{-1}y_n^{-1}+y_1y_2^2\cdot \dots \cdot y_{n-1}^{-1}}{y_1-y_1^{-1}}
y_n(y_n-y_n^{-1}){\det}_n''=\\
=&\frac{-y_1^{-1}y_n^{-1}+y_1y_2^2\cdot \dots \cdot y_{n-1}^{-1}}{y_1-y_1^{-1}}
(y_n-y_n^{-1}){\det}_n'' \left(\frac{-y_{n-1}(y_{n-2}-y_{n-2}^{-1})}{y_{n-1}-y_{n-1}^{-1}}
{\det}_{n-2}''-\frac{y_{n-1}^{-1}-y_{n-1}y_n}{y_n-y_n^{-1}}{\det}_{n-1}'' \right).\ea \] Using
the induction hypothesis we get an expression in the five variables $y_1,y_2^2\cdot \dots \cdot
y_{n-3}^2,y_{n-2},y_{n-1},y_n$ which can be computed to equal $ y_1^2\cdot \dots \cdot y_n^2+
y_1^{-2}\cdot \dots \cdot y_n^{-2}-2$.

\subsection{Proof of Proposition \ref{mainprop}}
Let $\e_2,\dots,\e_l \in \{-1,+1\}$.
Let $\varphi:\Z[y_1^{\pm 1},\dots,y_l^{\pm 1}] \to
\Z[y_1^{\pm 1},\dots,y_l^{\pm 1}]$ be the ring homomorphism
induced by $\varphi(y_1)=y_1$
and $\varphi(y_i)=y_i^{\e_i}, i=2,\dots,l$, denote the induced map
on $M_{2m \times 2m}(\Z[y_1^{\pm 1},\dots,y_l^{\pm 1}])$
by $\varphi$ as well.
Write $B(\e_2,\dots,\e_l)$
for $B(v_1,\dots,v_l,\e_2,\dots,\e_l)$.

We see that if we multiply the $(2i-1)$--st and the $2i$--th row of
$Y_lB(\eps_2,\dots,\eps_l)-Y_l^{-1}B(\eps_2,\dots,\eps_l)^t$ by $\eps_i$, $i=2,\dots,l$, then we
get $\varphi(Y_lB(1,\dots,1)-Y_l^{-1}B(1,\dots,1)^t)$, in particular the determinants are the
same, i.e.
\[ \ba{rl}
&\det(Y_lB(\eps_2,\dots,\eps_l)-Y_l^{-1}B(\eps_2,\dots,\eps_l)^t)\\
=&\varphi(\det(Y_lB(1,\dots,1)-Y_l^{-1}B(1,\dots,1)^t))\\
=&\varphi\left(1-2v_1+ \sum_{j=1}^l (v_j-v_{j+1})
\left(\prod_{i=1}^j y_i^2+\prod_{i=1}^j y_i^{-2}\right)\right)\\
=&1-2v_1+ \sum_{j=1}^l (v_j-v_{j+1}) \left(y_1^2\prod_{i=2}^j y_i^{2\e_i}+y_1^{-2}\prod_{i=2}^j
y_i^{-2\e_i}\right).
\ea
 \]
This proves Proposition \ref{mainprop}.

\end{document}